\documentclass[12pt]{article}

\usepackage{amsfonts,amssymb,indent,subeqnarray}
\usepackage{bm}    
\usepackage{url}   

\date{January 30, 2009 \\[1.5mm] revised April 29, 2009}

\begin{document}

\title{Real-variables characterization of \\
       generalized Stieltjes functions}

\author{
     {\small Alan D.~Sokal\thanks{Also at Department of Mathematics,
           University College London, London WC1E 6BT, England.}}  \\[-2mm]
     {\small\it Department of Physics}       \\[-2mm]
     {\small\it New York University}         \\[-2mm]
     {\small\it 4 Washington Place}          \\[-2mm]
     {\small\it New York, NY 10003 USA}      \\[-2mm]
     {\small\tt sokal@nyu.edu}               \\[-2mm]
     {\protect\makebox[5in]{\quad}}  
     \\[-2mm]
}

\maketitle
\thispagestyle{empty}   

\vspace*{-5mm}

\begin{abstract}
We obtain a characterization of generalized Stieltjes functions
of any order~$\lambda > 0$
in terms of inequalities for their derivatives on $(0,\infty)$.
When $\lambda=1$, this provides a new and simple proof
of a characterization of Stieltjes functions
first obtained by Widder in 1938.
\end{abstract}


\bigskip
\noindent
{\bf Mathematics Subject Classification (MSC 2000) codes:}
26A48 (Primary);
26D10, 30E05, 30E20, 44A10, 44A15, 44A60 (Secondary).

\clearpage

\newtheorem{theorem}{Theorem}

\renewcommand{\theenumi}{\alph{enumi}}
\renewcommand{\labelenumi}{(\theenumi)}
\def\eop{\hbox{\kern1pt\vrule height6pt width4pt
depth1pt\kern1pt}\medskip}
\def\prf{\par\noindent{\bf Proof.\enspace}\rm}
\def\rmk{\par\medskip\noindent{\bf Remark\enspace}\rm}

\newcommand{\be}{\begin{equation}}
\newcommand{\ee}{\end{equation}}
\newcommand{\<}{\langle}
\renewcommand{\>}{\rangle}
\newcommand{\widebar}{\overline}
\def\reff#1{(\protect\ref{#1})}
\def\spose#1{\hbox to 0pt{#1\hss}}
\def\ltapprox{\mathrel{\spose{\lower 3pt\hbox{$\mathchar"218$}}
    \raise 2.0pt\hbox{$\mathchar"13C$}}}
\def\gtapprox{\mathrel{\spose{\lower 3pt\hbox{$\mathchar"218$}}
    \raise 2.0pt\hbox{$\mathchar"13E$}}}
\def\textprime{${}^\prime$}
\def\proof{\par\medskip\noindent{\sc Proof.\ }}
\def\firstproof{\par\medskip\noindent{\sc First Proof.\ }}
\def\secondproof{\par\medskip\noindent{\sc Second Proof.\ }}
\def\qed{ $\square$ \bigskip}
\def\proofof#1{\bigskip\noindent{\sc Proof of #1.\ }}
\def\firstproofof#1{\bigskip\noindent{\sc First Proof of #1.\ }}
\def\secondproofof#1{\bigskip\noindent{\sc Second Proof of #1.\ }}
\def\altproofof#1{\bigskip\noindent{\sc Alternate Proof of #1.\ }}
\def\half{ {1 \over 2} }
\def\third{ {1 \over 3} }
\def\twothird{ {2 \over 3} }
\def\smfrac#1#2{{\textstyle{#1\over #2}}}
\def\smhalf{ \smfrac{1}{2} }
\newcommand{\real}{\mathop{\rm Re}\nolimits}
\renewcommand{\Re}{\mathop{\rm Re}\nolimits}
\newcommand{\imag}{\mathop{\rm Im}\nolimits}
\renewcommand{\Im}{\mathop{\rm Im}\nolimits}
\newcommand{\sgn}{\mathop{\rm sgn}\nolimits}
\newcommand{\tr}{\mathop{\rm tr}\nolimits}
\newcommand{\per}{\mathop{\rm per}\nolimits}
\newcommand{\supp}{\mathop{\rm supp}\nolimits}
\newcommand{\diag}{\mathop{\rm diag}\nolimits}
\def\hboxscript#1{ {\hbox{\scriptsize\em #1}} }
\renewcommand{\emptyset}{\varnothing}

\newcommand{\restrict}{\upharpoonright}
\newcommand{\implies}{\Longrightarrow}
\renewcommand{\iff}{\Longleftrightarrow}

\newcommand{\scra}{{\mathcal{A}}}
\newcommand{\scrb}{{\mathcal{B}}}
\newcommand{\scrc}{{\mathcal{C}}}
\newcommand{\scrd}{{\mathcal{D}}}
\newcommand{\scrf}{{\mathcal{F}}}
\newcommand{\scrg}{{\mathcal{G}}}
\newcommand{\scrh}{{\mathcal{H}}}
\newcommand{\scrk}{{\mathcal{K}}}
\newcommand{\scrl}{{\mathcal{L}}}
\newcommand{\scrm}{{\mathcal{M}}}
\newcommand{\scro}{{\mathcal{O}}}
\newcommand{\scrp}{{\mathcal{P}}}
\newcommand{\scrr}{{\mathcal{R}}}
\newcommand{\scrs}{{\mathcal{S}}}
\newcommand{\scrt}{{\mathcal{T}}}
\newcommand{\scrv}{{\mathcal{V}}}
\newcommand{\scrw}{{\mathcal{W}}}
\newcommand{\scrz}{{\mathcal{Z}}}

\newcommand{\am}{{\mathcal{AM}}}
\newcommand{\cm}{{\mathcal{CM}}}

\newcommand{\ahat}{{\widehat{a}}}
\newcommand{\Zhat}{{\widehat{Z}}}
\newcommand{\cc}{{\mathbf{c}}}
\renewcommand{\k}{{\mathbf{k}}}
\newcommand{\m}{{\mathbf{m}}}
\newcommand{\n}{{\mathbf{n}}}
\newcommand{\vv}{{\mathbf{v}}}
\newcommand{\bv}{{\mathbf{v}}}
\newcommand{\w}{{\mathbf{w}}}
\newcommand{\x}{{\mathbf{x}}}
\newcommand{\y}{{\mathbf{y}}}
\newcommand{\z}{{\bf z}}
\newcommand{\zbar}{{\bar{\mathbf{z}}}}
\newcommand{\zero}{{\mathbf{0}}}
\newcommand{\one}{{\mathbf{1}}}
\newcommand{\bR}{{\bf R}}  
\newcommand{\bRtilde}{{\widetilde{\bf R}}}
\newcommand{\bRhat}{{\widehat{\bf R}}}

\newcommand{\ofo}{ {{}_1 \! F_1} }
\newcommand{\oft}{ {{}_1 \! F_2} }

\newcommand{\C}{{\mathbb C}}
\newcommand{\D}{{\mathbb D}}
\newcommand{\Z}{{\mathbb Z}}
\newcommand{\N}{{\mathbb N}}
\newcommand{\Q}{{\mathbb Q}}
\newcommand{\R}{{\mathbb R}}
\newcommand{\HH}{{\mathbb H}}
\newcommand{\RR}{{\mathbb R}}

\newcommand{\varphibar}{{\bar{\varphi}}}
\newcommand{\bvarphi}{{\boldsymbol{\varphi}}}
\newcommand{\bvarphibar}{{\bar{\boldsymbol{\varphi}}}}
\newcommand{\psibar}{{\bar{\psi}}}

\def\cbar{{\overline{C}}}

\newcommand{\bigdash}{%
\smallskip\begin{center} \rule{5cm}{0.1mm} \end{center}\smallskip}


\newenvironment{sarray}{
             \textfont0=\scriptfont0
             \scriptfont0=\scriptscriptfont0
             \textfont1=\scriptfont1
             \scriptfont1=\scriptscriptfont1
             \textfont2=\scriptfont2
             \scriptfont2=\scriptscriptfont2
             \textfont3=\scriptfont3
             \scriptfont3=\scriptscriptfont3
           \renewcommand{\arraystretch}{0.7}
           \begin{array}{l}}{\end{array}}

\newenvironment{scarray}{
             \textfont0=\scriptfont0
             \scriptfont0=\scriptscriptfont0
             \textfont1=\scriptfont1
             \scriptfont1=\scriptscriptfont1
             \textfont2=\scriptfont2
             \scriptfont2=\scriptscriptfont2
             \textfont3=\scriptfont3
             \scriptfont3=\scriptscriptfont3
           \renewcommand{\arraystretch}{0.7}
           \begin{array}{c}}{\end{array}}

%
%
\newcommand{\stirlingsubset}[2]{\genfrac{\{}{\}}{0pt}{}{#1}{#2}}
\newcommand{\stirlingcycle}[2]{\genfrac{[}{]}{0pt}{}{#1}{#2}}
\newcommand{\assocstirlingsubset}[3]{%
{\genfrac{\{}{\}}{0pt}{}{#1}{#2}}_{\! \ge #3}}
\newcommand{\assocstirlingcycle}[3]{{\genfrac{[}{]}{0pt}{}{#1}{#2}}_{\ge #3}}
\newcommand{\euler}[2]{\genfrac{\langle}{\rangle}{0pt}{}{#1}{#2}}
\newcommand{\eulergen}[3]{{\genfrac{\langle}{\rangle}{0pt}{}{#1}{#2}}_{\! #3}}
\newcommand{\eulersecond}[2]{\left\langle\!\! \euler{#1}{#2} \!\!\right\rangle}
\newcommand{\eulersecondgen}[3]{%
{\left\langle\!\! \euler{#1}{#2} \!\!\right\rangle}_{\! #3}}
\newcommand{\binomvert}[2]{\genfrac{\vert}{\vert}{0pt}{}{#1}{#2}}


A real-valued function $f$ defined on an open interval $I \subseteq \R$
is said to be {\em completely monotone}\/ if it is $C^\infty$
and satisfies $(-1)^n f^{(n)}(x) \ge 0$ for all $x \in I$ and all $n \ge 0$.
The most important case is $I = (0,\infty)$,
where the Bernstein--Hausdorff--Widder theorem
\cite{Bernstein_29,Hausdorff_21a,Hausdorff_21b,Widder_31,Widder_46}
states that $f$ is completely monotone on $(0,\infty)$ if and only if
it can be written as the Laplace transform of a nonnegative measure
supported on $[0,\infty)$, i.e.
\be
   f(x)  \;=\;  \int\limits_{[0,\infty)} \! e^{-tx} \, d\mu(t)
 \label{eq.BHW}
\ee
with $\mu \ge 0$ and the integral convergent for all $x > 0$.\footnote{
   The book of Widder \cite{Widder_46}
   gives several different proofs of the
   Bernstein--Hausdorff--Widder theorem:
   one based on the Hausdorff moment problem
   and Carlson's theorem on analytic functions (pp.~160--161);
   one based on the Hausdorff moment problem and its uniqueness (pp.~162--163);
   one based on Laguerre polynomials (pp.~168--177);
   and one based on a real inversion formula for the Laplace transform
   (pp.~310--312).
   See also \cite[Chapter~I]{Donoghue_74} for a proof based on
   Newtonian interpolation polynomials,
   and \cite{Choquet_69} \cite[Chapter~2]{Phelps_01}
   for beautiful proofs based on Choquet theory.
}
Clearly, any such $f$ has an analytic continuation
to the right half-plane $\real x > 0$.\footnote{
   This latter property of completely monotone functions
   is fairly easy to prove by a direct argument
   not needing the full Bernstein--Hausdorff--Widder theorem:
   see e.g.\ \cite[pp.~146--147]{Widder_46}
   or \cite[pp.~13--14]{Donoghue_74}.
}

A real-valued function $f$ defined on $(0,\infty)$
is said to be a {\em Stieltjes function}\/ \cite{Stieltjes_1894}\footnote{
   More information on Stieltjes functions can be found in
   \cite[pp.~126--128]{Akhiezer_65} \cite{Berg_79,Berg_80}
   and the references cited therein.
}
if it can be written as a nonnegative constant plus
the Stieltjes transform \cite{Widder_38,Widder_46}
of a nonnegative measure supported on $[0,\infty)$, i.e.
\be
   f(x)  \;=\;  C \,+ \int\limits_{[0,\infty)} \! {d\rho(t) \over x+t}
\ee
with $C \ge 0$, $\rho \ge 0$
and the integral convergent for some (hence all) $x > 0$.
Clearly, every Stieltjes function is completely monotone on $(0,\infty)$,
but not every completely monotone function is Stieltjes.
It is thus of interest to obtain a characterization of
Stieltjes functions in terms of inequalities for the
derivatives of $f$ on $(0,\infty)$,
analogous to but stronger than
the inequalities defining complete monotonicity.
Such a characterization was obtained by Widder \cite{Widder_38} in 1938
(see also \cite[Chapter~VIII]{Widder_46}),
who proved (here $D=d/dx$):

\begin{theorem}
 \label{thm1.1}
  {$\!$ \bf \protect\cite[Theorem~12.5 and Lemma~12.52]{Widder_38} \ }
Let $f$ be a real-valued function defined on $(0,\infty)$.
Then the following are equivalent:
\begin{itemize}
   \item[(a)]  $f$ is a Stieltjes function.
   \item[(b)]  $f$ is $C^\infty$, and the quantities
\begin{subeqnarray}
   F_{n,k}(x)
   & = &
   (-1)^n \sum_{j=0}^k {k \choose j} \,
                       {(n+k)! \over (n+j)!} \:
                   x^j \, f^{(n+j)}(x)
       \slabel{def.Ftildenk.sum}  \\[2mm]
   & = &
      (-1)^n \, x^{-n} D^k x^{n+k} D^n f(x)
       \slabel{def.Ftildenk.deriv1}  \\[2mm]
   & = &
      (-1)^n \, D^{n+k} x^k f(x)
       \slabel{def.Ftildenk.deriv2}
       \label{def.Ftildenk}
\end{subeqnarray}
   are nonnegative for all $n,k \ge 0$ and all $x > 0$.
   \item[(c)]  $f$ is $C^\infty$, and we have $F_{0,0}(x) \ge 0$
      and $F_{k-1,k}(x) \ge 0$ for all $k \ge 1$ and all $x > 0$.
\end{itemize}
\end{theorem}

Since $F_{n,0} = (-1)^n f^{(n)}$,
the condition (b) is manifestly a strengthening of complete monotonicity.
The equivalence of the three formulae for $F_{n,k}$
is a straightforward computation.

{}From \reff{def.Ftildenk.deriv2} we see that
the nonnegativity of $F_{n,k}$ for all $n,k \ge 0$
is equivalent to the assertion that
all the functions $F_{0,k} = D^k x^k f$
are completely monotone on $(0,\infty)$.

It is fairly easy to see that (a) $\implies$ (b),
while (b) $\implies$ (c) is trivial.
Widder's proof of (c) $\implies$ (a) was, by contrast, fairly long,
and was based on explicit construction of a differential operator $L_{k,t}$
that provides a real inversion formula for the Stieltjes transform.
Along the way he also gave \cite[Lemma~12.52]{Widder_38}
a direct real-variables proof of (c) $\implies$ (b),
but he used this only for technical purposes,
to guarantee the complete monotonicity and hence the real-analyticity
of $f$ on $(0,\infty)$ \cite[p.~48]{Widder_38}.\footnote{
   See \cite[Chapter~VIII]{Widder_46} for a slightly different proof
   that does not make use of \cite[Lemma~12.52]{Widder_38}.
}

In addition, Widder \cite[Theorem~10.1]{Widder_36} proved, two years earlier,
a slight variant of Theorem~\ref{thm1.1}(a) $\iff$ (b)
--- treating the case in which the measure $\mu$ is required to be finite ---
by applying the Bernstein--Hausdorff--Widder theorem to the functions
$F_{0,k}$ and then analyzing the relationship between
the representing measures $\mu_k$.

In this paper I would like to give an extremely short and simple proof
of Theorem~\ref{thm1.1},
which moreover extends to provide a new characterization
of the generalized Stieltjes functions of any order~$\lambda > 0$
(see Theorem~\ref{thm1.2} below).
The key idea is to use the well-known solubility conditions
for the Hausdorff moment problem to prove (b) $\implies$ (a);
we then rely on \cite[Lemma~12.52]{Widder_38} for (c) $\implies$ (b).
Let us recall that a sequence ${\bf c} = (c_n)_{n=0}^\infty$
is said to be a {\em Hausdorff moment sequence}\/
if there exists a finite nonnegative measure $\nu$ on $[0,1]$ such that
\be
   c_n  \;=\;  \int\limits_{[0,1]} t^n \, d\nu(t)
   \quad\hbox{for all } n \ge 0  \;,
\ee
and it is said to be {\em completely monotone}\/ if
\be
   (-1)^k (\Delta^k {\bf c})_n 
   \;\equiv\;
   \sum_{j=0}^k (-1)^j {k \choose j} c_{n+j}
   \;\ge\;  0
   \quad\hbox{for all } n,k \ge 0  \;.
\ee
Hausdorff \cite{Hausdorff_21a} proved in 1921
that a sequence ${\bf c} = (c_n)_{n=0}^\infty$
is a Hausdorff moment sequence if and only if it is completely monotone;
furthermore, the representing measure $\nu$ is unique.\footnote{
   See also \cite[pp.~8--9]{Shohat_43},
   \cite[pp.~60--61 and 100--109]{Widder_46}
   or \cite[pp.~74--76]{Akhiezer_65}.
   ``Only if'' is quite easy;  proving ``if'' takes more work.
}
This is obviously a discrete analogue of the
Bernstein--Hausdorff--Widder theorem.

Our method also handles, with no extra work,
the generalized Stieltjes transform in which the kernel $1/(x+t)$
is replaced by $1/(x+t)^\lambda$ for some exponent $\lambda > 0$
\cite[Section~8]{Widder_38}
\cite{Pollard_42,Sumner_49,Byrne_74,Love_80,Love_82}.
Let us say that a real-valued function $f$ on $(0,\infty)$
is a {\em generalized Stieltjes function of order~$\lambda$}\/
(and write $f \in \scrs_\lambda$)
if it can be written in the form
\be
   f(x)  \;=\;  C \,+ \int\limits_{[0,\infty)} \! {d\rho(t) \over (x+t)^\lambda}
 \label{def.Slambda}
\ee
with $C \ge 0$, $\rho \ge 0$
and the integral convergent for some (hence all) $x > 0$.
Since
\be
   {1 \over (x+t)^\lambda}
   \;=\;
   {\Gamma(\lambda') \over \Gamma(\lambda) \, \Gamma(\lambda'-\lambda)}
   \int\limits_0^\infty \!
    u^{\lambda'-\lambda-1} \,  {1 \over (x+t+u)^{\lambda'}} \: du
\ee
whenever $\lambda < \lambda'$,
it follows that $\scrs_\lambda \subseteq \scrs_{\lambda'}$
whenever $\lambda \le \lambda'$.
It is also suggestive that the representation \reff{def.Slambda} tends formally
as $\lambda \uparrow \infty$ to the representation \reff{eq.BHW}
characteristic of complete monotonicity,
in the sense that
$\lim\limits_{\lambda \uparrow \infty}
 (\lambda t)^\lambda / (x+ \lambda t)^\lambda = e^{-x/t}$.

We shall prove the following real-variables characterization
of the generalized Stieltjes functions of order~$\lambda$:

\begin{theorem}
 \label{thm1.2}
Let $\lambda > 0$,
and let $f$ be a real-valued function defined on $(0,\infty)$.
Then the following are equivalent:
\begin{itemize}
   \item[(a)]  $f$ is a generalized Stieltjes function of order~$\lambda$.
   \item[(b)]  $f$ is $C^\infty$, and the quantities
\begin{subeqnarray}
   F^{[\lambda]}_{n,k}(x)
   & = &
   (-1)^n \sum_{j=0}^k {k \choose j} \,
                       {\Gamma(n+k+\lambda) \over \Gamma(n+j+\lambda)} \:
                   x^j \, f^{(n+j)}(x)  \slabel{eq.thm1.2.a} \\[2mm]
   & = &
      (-1)^n \, x^{-(n+\lambda-1)} D^k x^{n+k+\lambda-1} D^n f(x)
                                        \slabel{eq.thm1.2.b}
\end{subeqnarray}
   are nonnegative for all $n,k \ge 0$ and all $x > 0$.
\end{itemize}
\end{theorem}

\noindent
When $\lambda=1$ this reduces to Theorem~\ref{thm1.1}(a,b).

Since $F^{[\lambda]}_{n,0} = (-1)^n f^{(n)}$,
the condition (b) is manifestly a strengthening of complete monotonicity.
Furthermore, $F^{[\lambda]}_{n,k}(x)$
is a polynomial in $\lambda$ of degree $k$,
with leading coefficient
\be
   \lim\limits_{\lambda\to\infty} 
   {F^{[\lambda]}_{n,k}(x)  \over  \lambda^k}
   \;=\;
   (-1)^n \, f^{(n)}(x)
   \;.
\ee
So condition (b) tends formally as $\lambda \uparrow \infty$
to the definition of complete monotonicity,
and Theorem~\ref{thm1.2} tends formally
to the Bernstein--Hausdorff--Widder theorem.
At the other extreme, we have
\begin{subeqnarray}
   \lim\limits_{\lambda\to 0} F^{[\lambda]}_{0,1}(x) 
   & = &  x f'(x)  \\[2mm]
   \lim\limits_{\lambda\to 0} F^{[\lambda]}_{1,0}(x) 
   & = &   -f'(x)
\end{subeqnarray}
so that the only functions that are generalized Stieltjes
of all orders $\lambda > 0$ are the nonnegative constants.

\bigskip

{\bf Remarks.}
1.  The equivalence of the two formulae for $F^{[\lambda]}_{n,k}$
in \reff{eq.thm1.2.a}/\reff{eq.thm1.2.b}
is a straightforward computation.
However, for $\lambda \neq 1$ we do not know any simple rewriting of
$F^{[\lambda]}_{n,k}(x)$
analogous to the third formula \reff{def.Ftildenk.deriv2},
nor do we know (except possibly for integer values of $\lambda$, see below)
any characterization of the generalized Stieltjes functions
in terms of a proper subset of the $\{F^{[\lambda]}_{n,k}\}$
analogous to Theorem~\ref{thm1.1}(c).
Even when $\lambda=1$, it is an interesting open question
to find other proper subsets of the $\{F^{[\lambda]}_{n,k}\}$,
besides the one given in Theorem~\ref{thm1.1}(c),
whose nonnegativity is equivalent to that of the whole set.

2. It would also be interesting to show directly
that the conditions (b) get weaker as $\lambda$ grows.
The most obvious approach would be to write all the derivatives
$(\partial^\ell/\partial \lambda^\ell) F^{[\lambda]}_{n,k}$
as nonnegative linear combinations of $\{F^{[\lambda]}_{n',k'}\}$.

3. Some of Widder's results \cite[Theorems~8.2 and 8.3]{Widder_38}
may imply an alternative characterization of the generalized Stieltjes functions
of order~$\lambda$ that generalizes that of Theorem~\ref{thm1.1}(c).
When $\lambda$ is an integer,
this characterization will apparently involve the condition that
$F^{[\lambda]}_{k-\lambda,k}(x) \ge 0$ for all $k \ge \lambda$ and all $x > 0$,
probably together with the nonnegativity of a few other $F^{[\lambda]}_{n,k}$
(e.g.\ $F^{[\lambda]}_{0,0}$).
When $\lambda$ is noninteger, however,
this characterization will be nonlocal,
involving convolution as well as differentiation.
\qed

\proofof{Theorem~\ref{thm1.2}}
(a) $\implies$ (b):
Suppose that
\be
   f(x)  \;=\;  C \,+ \int\limits_{[0,\infty)} \! {d\rho(t) \over (x+t)^\lambda}
 \label{eq.f.stieltjes}
\ee
with $C \ge 0$, $\rho \ge 0$ and
$\int \! d\rho(t)/(1+t)^\lambda < \infty$.
Then $f$ is infinitely differentiable on $(0,\infty)$, with
\be
   f^{(n)}(x)  \;=\;
   C \delta_{n,0} \,+\,
   (-1)^n \, {\Gamma(n+\lambda) \over \Gamma(\lambda)}
   \! \int\limits_{[0,\infty)} \!\! {d\rho(t) \over (x+t)^{n+\lambda}}
   \qquad\hbox{for all } n \ge 0 \;.
\ee
It follows that
\be
   f^{[\lambda]}_n(x)
   \;\equiv\;
   (-1)^n \, {\Gamma(\lambda) \over \Gamma(n+\lambda)} \: x^n \, f^{(n)}(x)
   \;=\;  \int\limits_{[0,1]} u^n \, d\nu_x(u)
   \;,
\ee
where $d\nu_x(u)$ is the image of the measure $d\rho(t)/(x+t)^\lambda$
under the map $u = (1+t/x)^{-1}$
together with a point mass $C$ at $u=0$.
In other words, for each $x>0$ the sequence
$\bm{f}^{[\lambda]}(x) = (f^{[\lambda]}_n(x))_{n=0}^\infty$
is a Hausdorff moment sequence;
therefore, by (the easy half of) Hausdorff's theorem,
the sequence $\bm{f}^{[\lambda]}(x)$ is completely monotone,
i.e. the functions
\be
   f^{[\lambda]}_{n,k}(x)
   \;\equiv\;
      (-1)^k \, [\Delta^k \bm{f}^{[\lambda]}(x)]_n
   \;=\;
      (-1)^n x^n \sum_{j=0}^k {k \choose j} \,
                   {\Gamma(\lambda) \over \Gamma(n+j+\lambda)} \:
                   x^j \, f^{(n+j)}(x)
\ee
are nonnegative for all $n,k \ge 0$ and all $x>0$.
The same is therefore true of the functions
\be
   F^{[\lambda]}_{n,k}(x)
   \;\equiv\;
   {\Gamma(n+k+\lambda) \over \Gamma(\lambda)}  \:
       {f^{[\lambda]}_{n,k}(x) \over x^n}
   \;\,.
\ee
This proves (a) $\implies$ (b).

(b) $\implies$ (a):
Now we use the sufficiency half of Hausdorff's theorem:
it follows that, for each $x>0$, there exists a finite nonnegative measure
$\nu_x$ on $[0,1]$ such that
\be
  (-1)^n \, {\Gamma(\lambda)  \over \Gamma(n+\lambda)}
  \:  x^n \, f^{(n)}(x)
  \;=\;  \int\limits_{[0,1]} u^n \, d\nu_x(u)
   \qquad\hbox{for all } n \ge 0 \;.
 \label{eq.proof.hausdorff}
\ee
Changing variables back to $t = x(u^{-1} - 1)$,
we see that there exists a nonnegative measure
$\rho_x$ on $[0,\infty)$
satisfying $\int \! d\rho_x(t)/(x+t)^\lambda < \infty$,
and a constant $C_x \ge 0$, such that
\be
   f^{(n)}(x)  \;=\;
   C_x \delta_{n,0} \,+\,
   (-1)^n \, {\Gamma(n+\lambda) \over \Gamma(\lambda)}
   \! \int\limits_{[0,\infty)} \!\! {d\rho_x(t) \over (x+t)^{n+\lambda}}
   \qquad\hbox{for all } n \ge 0
 \label{eq.star1}
\ee
[namely, $d\rho_x(t) = (x+t)^\lambda d\varphi_x(\nu_x)(t)$
where $\varphi_x(u) = x(u^{-1} - 1)$,
and $C_x = \nu_x(\{0\})$].
We now use the fact that (b) implies the complete monotonicity of $f$,
hence the existence of an analytic continuation of $f$
to the right half-plane;
in particular, the Taylor series for $f$ or any of its derivatives
around the point $x$ must have radius of convergence at least $x$.
So let us take \reff{eq.star1} with $n$ replaced by $n+k$,
multiply it by $\xi^k/k!$, and sum over $k \ge 0$:
for $|\xi| < x$ the series is absolutely convergent, and we obtain
\be
   f^{(n)}(x+\xi)
   \;=\;
   C_x \delta_{n,0} \,+\,
   (-1)^n \, {\Gamma(n+\lambda) \over \Gamma(\lambda)}
   \! \int\limits_{[0,\infty)}  \!
      {d\rho_x(t) \over (x+\xi+t)^{n+\lambda}}
   \qquad\hbox{for all } n \ge 0
\ee
whenever $\xi \in (-x,x)$,
or in other words
\be
   f^{(n)}(y)
   \;=\;
   C_x \delta_{n,0} \,+\,
   (-1)^n \, {\Gamma(n+\lambda) \over \Gamma(\lambda)}
   \! \int\limits_{[0,\infty)}  \! {d\rho_x(t) \over (y+t)^{n+\lambda}}
   \qquad\hbox{for all } n \ge 0
\ee
whenever $y \in (0,2x)$,
or equivalently
\be
  (-1)^n \, {\Gamma(\lambda)  \over \Gamma(n+\lambda)}
  \:  y^n \, f^{(n)}(y)
  \;=\;  \int\limits_{[0,1]} u^n \, d\nu'_{x,y}(u)
   \qquad\hbox{for all } n \ge 0
 \label{eq.proof.1}
\ee
where $d\nu'_{x,y}(u)$ is the image of the measure $d\rho_x(t)/(y+t)^\lambda$
under the map $u = (1+t/y)^{-1}$
together with a point mass $C_x$ at $u=0$.
On the other hand, we already know from \reff{eq.proof.hausdorff} that
\be
  (-1)^n \, {\Gamma(\lambda)  \over \Gamma(n+\lambda)}
  \:  y^n \, f^{(n)}(y)
  \;=\;  \int\limits_{[0,1]} u^n \, d\nu_y(u)
   \qquad\hbox{for all } n \ge 0  \;.
 \label{eq.proof.2}
\ee
Comparing \reff{eq.proof.1}/\reff{eq.proof.2},
we see that the measures $\nu'_{x,y}$ and $\nu_y$ have the same moments
whenever $0 < y < 2x$;
so by the uniqueness in the Hausdorff moment problem,
we conclude that $\nu'_{x,y} = \nu_y$
and hence $C_x = C_y$ and $\rho_x = \rho_y$ whenever $0 < y < 2x$.
In particular, $C_x = C_y$ and $\rho_x = \rho_y$ whenever $0 < y < x$;
and this implies, using the symmetry $x \leftrightarrow y$,
that $C_x = C_y$ and $\rho_x = \rho_y$ for all $x,y > 0$.
This proves (b) $\implies$ (a).
\qed

{\bf Remark.}
Here is an alternate proof of (a) $\implies$ (b):
since \cite[p.~299]{Love_80}
\be
   (-1)^n \, x^{-(n+\lambda-1)} \, {d^k \over dx^k} \, x^{n+k+\lambda-1} \,
    {d^n \over dx^n} \, {1 \over (x+t)^\lambda}
   \;=\;
   {\Gamma(n+k+\lambda) \over \Gamma(\lambda)} \:
       {t^k \over (x+t)^{n+k+\lambda}}
   \;,
\ee
the representation \reff{eq.f.stieltjes} implies that
\be
   F^{[\lambda]}_{n,k}(x)
   \;=\;
   {\Gamma(n+k+\lambda) \over \Gamma(\lambda)}  \:
   \Biggl[ C \delta_{n,0} \;+
      \int\limits_{[0,\infty)} \! {t^k \over (x+t)^{n+k+\lambda}} \, d\rho(t)
   \Biggr]
   \;\ge\; 0  \;.
\vspace*{-6mm}
\ee
\qed

\bigskip

Let us conclude by remarking that the Stieltjes functions
also have a beautiful complex-analysis characterization:
a function $f \colon\, (0,\infty) \to \R$
is Stieltjes if and only if it is the restriction to $(0,\infty)$
of an analytic function on the cut plane $\C \setminus (-\infty,0]$
satisfying $f(z) \ge 0$ for $z>0$ and
$\Im f(z) \le 0$ for $\Im z > 0$.
See e.g.\ \cite[p.~127]{Akhiezer_65} or \cite{Berg_80}.
It would be interesting to know whether the generalized Stieltjes functions
of order~$\lambda$ have an analogous complex-analysis characterization
for some (or all) $\lambda \neq 1$.

\section*{Acknowledgments}

I wish to thank Rob Corless for drawing my attention to Stieltjes functions,
and especially Christian Berg for extremely helpful correspondence.
I also wish to thank an anonymous referee for helpful suggestions
concerning the exposition.

Finally, I thank the Institut Henri Poincar\'e -- Centre Emile Borel
for hospitality during the programme on
Interacting Particle Systems, Statistical Mechanics and Probability Theory
(September--December 2008),
where most of this work was carried out.

This research was supported in part by
U.S.\ National Science Foundation grant PHY--0424082.

\end{document}